\documentstyle[12pt]{article}
\begin{document}
\begin{titlepage}
\begin{flushright}
q-alg/9712011
\end{flushright}
\vskip.3in

\begin{center}
{\Large \bf On Super RS Algebra and Drinfeld Realization of Quantum
Affine Superalgebras} 
\vskip.3in
{\large Mark D. Gould and Yao-Zhong Zhang}
\vskip.2in
{\em Department of Mathematics, University of Queensland, Brisbane,
     Qld 4072, Australia

Email: yzz@maths.uq.edu.au}
\end{center}

\vskip 2cm
\begin{center}
{\bf Abstract}
\end{center}
We describe the realization of the super Reshetikhin-Semenov-Tian-Shansky
(RS) algebra in quantum affine superalgebras, thus generalizing the 
approach of Frenkel-Reshetikhin to the supersymmetric (and twisted) case. 
The algebraic homomorphism
between the super RS algebra and the Drinfeld current realization
of quantum affine superalgebras is established by using the Gauss 
decomposition technique of Ding-Frenkel. As an application, we obtain
Drinfeld realization of quantum affine superalgebra $U_q[osp(1|2)^{(1)}]$ 
and its degeneration -- central extended super Yangian double
$DY_\hbar[osp(1|2)^{(1)}]$.

\vskip 3cm
\noindent{\bf Mathematics Subject Classifications (1991):} 81R10, 17B37, 16W30

\end{titlepage}

%  Greek letters

\def\a{\alpha}
\def\b{\beta}
\def\d{\delta}
\def\e{\epsilon}
\def\g{\gamma}
\def\k{\kappa}
\def\l{\lambda}
\def\o{\omega}
\def\t{\theta}
\def\s{\sigma}
\def\D{\Delta}
\def\L{\Lambda}

\def\G{{\cal G}}
\def\Gk{{\cal G}^{(k)}}
\def\R{{\cal R}}
\def\hR{{\hat{\cal R}}}
\def\C{{\bf C}}
\def\P{{\bf P}}

\def\uqgh{{U_q[gl(m|n)^{(1)}]}}
\def\uqoh{{U_q[osp(1|2)^{(1)}]}}

% Shorthands for \begin{equation} and the like

\def\beq{\begin{equation}}
\def\eeq{\end{equation}}
\def\bea{\begin{eqnarray}}
\def\eea{\end{eqnarray}}
\def\ba{\begin{array}}
\def\ea{\end{array}}
\def\no{\nonumber}
\def\lt{\left}
\def\rt{\right}
\newcommand{\bq}{\begin{quote}}
\newcommand{\eq}{\end{quote}}

\newtheorem{Theorem}{Theorem}
\newtheorem{Definition}{Definition}
\newtheorem{Proposition}{Proposition}
\newtheorem{Lemma}[Theorem]{Lemma}
\newtheorem{Corollary}[Theorem]{Corollary}
\newcommand{\proof}[1]{{\bf Proof. }
        #1\begin{flushright}$\Box$\end{flushright}}

\newcommand{\sect}[1]{\setcounter{equation}{0}\section{#1}}
\renewcommand{\theequation}{\thesection.\arabic{equation}}

\sect{Introduction\label{intro}}

Recently we introduced \cite{Zha97a}, without elaboration,
a super version of the
Reshetikhin-Semenov-Tian-Shansky (RS) algebra
\cite{Res90}. Using this super RS algebra and
a super analogue of the Gauss decomposition formula of Ding-Frenkel 
\cite{Din93}, we obtain Drinfeld realization \cite{Dri88} of quantum affine
superalgebra $\uqgh$ 
(see also \cite{Cai97} for the special case of $m=n=1$). 

In this paper, we adopt a more conceptual way of proceeding. We
shall show that the super RS algebra can be realized in quantum affine
superalgebras, which enables extension of the work by Frenkel-Reshetikhin
\cite{Fre92} to the supersymmetric (and twisted) case. The algebra homomorphism
between the super RS algebra and the Drinfeld realization of the
quantum affine superalgebras is achieved by using the Gauss
decomposition technique of Ding-Frenkel. As an application, we obtain
the Drinfeld realization of $\uqoh$ and its degenerated algebra --
central extended super-Yangian double $DY_\hbar[osp(1|2)^{(1)}]$.

\sect{Super RS Algebra and its Realization in Quantum Affine
      Superalgebras}

Let us start with the definition of the super RS algebra. Let 
$R(z)\in End(V\otimes V)$, where $V$ is a ${\bf Z}_2$ graded vector space,
be a matrix obeying the weight conservation condition
$R(z)_{\a\b,\a'\b'}\neq 0$  only when
$[\a']+[\b']+[\a]+[\b]=0$ mod$2$, and the
graded Yang-Baxter equation (YBE) 
\beq\label{rrr}
R_{12}(z)R_{13}(zw)R_{23}(w)=R_{23}(w)R_{13}(zw)R_{12}(z).
\eeq
The multiplication rule for the tensor product is defined for
homogeneous elements $a,~ b,~ a'$, $b'$ by
\beq
(a\otimes b)(a'\otimes b')=(-1)^{[b][a']}\,(aa'\otimes bb'),
\eeq
where $[a]\in{\bf Z}_2$ 
denotes the grading of the element $a$. We introduce \cite{Zha97a} 
\begin{Definition}\label{rs}: Let $R(\frac{z}{w})$ be a R-matrix
satisfying the graded YBE (\ref{rrr}). The super RS algebra $U(\R)$ is
generated by invertible $L^\pm(z)$, satisfying
\bea
R({z\over w})L_1^\pm(z)L_2^\pm(w)&=&L_2^\pm(w)L_1^\pm(z)R({z\over
         w}),\no\\
R({z_+\over w_-})L_1^+(z)L_2^-(w)&=&L_2^-(w)L_1^+(z)R({z_-\over
         w_+}),\label{super-rs}
\eea
where $L_1^\pm(z)=L^\pm(z)\otimes 1$, $L_2^\pm(z)=1\otimes L^\pm(z)$
and $z_\pm=zq^{\pm {c\over 2}}$. For the first formula of
(\ref{super-rs}), the expansion direction of $R({z\over w})$ can be
chosen in $z\over w$ or $w\over z$, but for the second formula, the
expansion direction must only be in $z\over w$.
\end{Definition}

The algebra $U(\R)$ is a graded Hopf algebra: its coproduct is defined by
\beq
\D(L^\pm(z)=L^\pm(zq^{\pm 1\otimes {c\over 2}})\stackrel{.}{\otimes}
     L^\pm(zq^{\mp {c\over 2}\otimes 1}),
\eeq
and its antipode is
\beq
S(L^\pm(z))=L^\pm(z)^{-1}.
\eeq

We now consider a realization of $U(\R)$ in quantum affine superalgebra
$U_q[\Gk]$, generalizing the Frenkel-Reshetikhin description to the
supersymmetric case.
Let us first of all recall some facts about the affine superalgebra
${\cal G}^{(k)}$. For simplicity, we restrict ourselves to the case of 
$k\leq 2$. Let ${\cal G}_0$ be the fixed point subalgebra under
the diagram automorphism $\hat{\tau}$ of ${\cal G}$ of order $k$. 
In the case of $k=1$, we have ${\cal G}_0\equiv{\cal G}$.
For $k=2$ we may decompose ${\cal G}$ as ${\cal G}_0$ plus a ${\cal
G}_0$-representation ${\cal G}_1$ of ${\cal G}$. Let
\beq
\psi=\lt\{
\begin{array}{l}
{\rm highest~ root~ of}~ {\cal G}_0\equiv {\cal G}, ~~{\rm for}~ k=1,\\
{\rm highest~ weight~ of~ the}~ {\cal G}_0-{\rm representation}~ 
  {\cal G}_1, ~~{\rm for}~ k=2.
\end{array}
\rt.
\eeq
Following the usual convention, we denote the weight
of ${\cal G}^{(k)}$ by $\L\equiv (\l,\kappa,\tau)$, where $\l$ is a
weight of ${\cal G}_0$. With this 
convention the nondegenerate form $(~,~)$ induced on the weights
can be expressed as
\beq
(\L,\L')=(\l,\l')+\kappa\tau'+\kappa'\tau.
\eeq
Let $h_{\hat{\rho}}$ denote the unique element of the Cartan subalgebra
of ${\cal G}^{(k)}$ satisfying $h_{\hat{\rho}}(\a_i)=\frac{1}{2}
(\a_i,\a_i)$, where
$\a_i~(0\leq i\leq r)$ are simple roots of ${\cal G}^{(k)}$. Then 
$h_{\hat{\rho}}$ is given by
\beq
h_{\hat{\rho}}=h_\rho+gd,
\eeq
where $g=\frac{k}{2}(\psi,\psi+2\rho)$,  $\rho$ is the graded half-sum of 
positive roots of ${\cal G}_0$ and $d$ is the usual level operator. 

We shall not give the relations obeyed by the simple 
generators $\{h_i,~E_i,~F_i,~d,~0\leq i\leq r\}$ of
$U_q[{\cal G}^{(k)}]$, 
but mention that $U_q[\Gk]$ is endowed with a graded Hopf algebra
structure with coproduct and antipode given by
\bea
\D(h_i)&=&h_i\otimes 1+1\otimes h_i,~~~~\D(d)=d\otimes 1+1\otimes
          d,\no\\
\D(E_i)&=&E_i\otimes q^{-\frac{h_i}{2}}+q^{\frac{h_i}{2}}\otimes
          E_i,~~~~
\D(F_i)=F_i\otimes q^{-\frac{h_i}{2}}+q^{\frac{h_i}{2}}\otimes
          F_i,\no\\
S(a)&=&-q^{-h_{\hat{\rho}}}\;a\;q^{h_{\hat{\rho}}},~~~~\forall 
       a=d,\;h_i,\;E_i,\;F_i.
\eea

We denote by $\hat{\R}$ the universal R-matrix of $U_q[\Gk]$, which 
satisfies the graded YBE:
\beq
\hR_{12}\hR_{13}\hR_{23}=\hR_{23}\hR_{13}\hR_{12}
\eeq
and the coproduct properties:
\bea
&&(\D\otimes 1)\hR=
  \hR_{13}\hR_{23},~~~~~(1\otimes\D)\hR=\hR_{13}\hR_{12},\no\\
&&(\D^T\otimes 1)\hR=
  \hR_{23}\hR_{13},~~~~~(1\otimes\D^T)\hR=\hR_{12}\hR_{13},\label{D-hR}
\eea
where $\D^T\equiv T\D$, $T$ is the graded twist map so that for $a=\sum_i
\,a_i,~b=\sum_j\,b_j$, 
\beq
T(a\otimes b)=\sum_{i,j}\,(-1)^{[a_i][b_j]}\;(b_j
  \otimes a_i).
\eeq

Let $\hR_{21}\equiv \hR^T\equiv T\hR$. Note that
$\lt(\hR^T\rt)^{-1}=\lt(\hR^{-1}\rt)^T$ also satisfies the coproduct
properties (\ref{D-R}). Both $\hR$ and $\lt(\hR^T\rt)^{-1}$ satisfy the
intertwining property
\beq
\hR\D(a)=\D^T(a)\hR,~~~~~\lt(\hR^T\rt)^{-1}\D(a)=\D^T(a)\lt(\hR^T\rt)^{-1},~~~~~
  \forall a\in U_q[\Gk].
\eeq

The Hopf superalgebra $U_q[\Gk]$ contains two important Hopf subalgebras
$U_q^+$ and $U_q^-$ which are generated by
$\{E_i,\;h_i,\;d|i=0,1,\cdots,r\}$ and
$\{F_i,\;h_i,\;d|i=0,1,\cdots,r\}$, respectively. By Drinfeld's 
quantum double construction, the
universal R-matrix  $\hR$ can be written in the form 
\beq
\hR=\lt(I\otimes I+\sum_t\,a^t\otimes a_t\rt)\cdot q^{-\sum_{i=1}^r
  H^i\otimes H_i-c\otimes d-d\otimes c},
\eeq
where $\{a^t\}\in U_q^+,~\{a_t\}\in U_q^-$ do not depend on $d$, and are 
generated by $\{E_i,\;h_i\}$ and $\{F_i,\;h_i\}$,
respectively; $c$ is given by $\frac{c}{k}=h_0+h_\psi$, and $\{H^i\},~
\{H_i\}~(i=1,2,\cdots,r)$ satisfy
\beq
\sum_{i=1}^r\L(H^i)\L'(H_i)=(\l,\l').
\eeq
Let us remark that another form of
the universal R-matrix corresponding to the same coproduct reads
\beq
\hR'=\lt(I\otimes I+\sum_t\,a'_t\otimes a'^t\rt)\cdot q^{\sum_{i=1}^r
  H_i\otimes H^i+c\otimes d+d\otimes c},
\eeq
where $\{a'_t\}\in U_q^-$ and $\{a'^t\}\in U_q^+$ do not depend on
$d$. $\hR'$ can be identified with
$\lt(\hR^T\rt)^{-1}$.

Following Frenkel and Reshetikhin \cite{Fre92}, we define the
`normalized' universal R-matrix by the formula
\beq
\R=\hR\cdot q^{c\otimes d+d\otimes c}.\label{normalized-R}
\eeq
Then it can be shown that $\R$ satisfies the following `normalized
YBE'
\beq
\R_{12}q^{-d\otimes c\otimes 1}\R_{13}q^{d\otimes c\otimes 1}\R_{23}
  =\R_{23}q^{-1\otimes c\otimes d}\R_{13}q^{1\otimes c\otimes d}
   \R_{12},\label{normalized-ybe1}
\eeq
and the coproduct properties:
\bea
&&(\D\otimes 1)\R=
  \R_{13}q^{-c\otimes 1\otimes d}\R_{23}q^{c\otimes 1\otimes d},~~~~~
  (1\otimes\D)\R=\R_{13}q^{-d\otimes 1\otimes c}\R_{12}q^{d\otimes
  1\otimes c},\no\\
&&(\D^T\otimes 1)\R=
  \R_{23}q^{-1\otimes c\otimes d}\R_{13}q^{1\otimes c\otimes d},~~~~~
  (1\otimes\D^T)\R=\R_{12}q^{-d\otimes c\otimes 1}\R_{13}q^{d\otimes
  c\otimes 1}.\label{D-R}
\eea

Let $z$ be a formal variable. 
Define an automorphism $D_z$ of $U_q[{\cal G}^{(k)}]$ by
\beq
D_z(e_i)=z^{\frac{\d_{i0}}{k}}e_i,~~~~~D_z(f_i)=z^{-\frac{\d_{i0}}{k}}f_i,~~~~
D_z(h_i)=h_i,~~~~D_z(d)=d.
\eeq
It is worth noting that the automorphism $D_z$ is related to the
generator $d$  by the formula:
\beq
D_z(a)=z^d\,a\,z^{-d},~~~~~~\forall a\in U_q[{\cal G}^{(k)}]
\eeq
which is easily checked by means of the commutation relations,
\beq
[d,e_i]=\frac{\d_{i0}}{k}e_i,~~~~~[d,f_i]=-\frac{\d_{i0}}{k}f_i,~~~~~
[d,h_i]=0.
\eeq
We define a universal R-matrix $\R(z)$ depending on the formal parameter
$z$ by the formula
\beq
\R(z)=(D_z\otimes 1)\R=(1\otimes D_{z^{-1}})\R.
\eeq
Then (\ref{normalized-ybe1}) implies the following relation for  $\R(z)$:
\beq
\R_{12}(z)\R_{13}(zwq^{-c_2})\R_{23}(w)=\R_{23}(w)\R_{13}(zwq^{c_2})
   \R_{12}(z),\label{normalized-ybe2}
\eeq
where $c_2=1\otimes c\otimes 1$.
One can also show that $\R(z)$  enjoys
\beq
(S\otimes 1)(\R(z))=\R(zq^{-c\otimes 1})^{-1},~~~~
(1\otimes S^{-1})(\R(z))=\R(zq^{1\otimes c})^{-1},\label{s-r}
\eeq
and $(S\otimes S)(\R(z))=\R(zq^{1\otimes c-c\otimes 1})$.

For a finite dimensional representation $\pi_V$ supplied by the graded
vector space $V$, we define $R(z)\in End(V\otimes V)$ by
\beq
R(z)=(\pi_V\otimes\pi_V)\R(z).
\eeq
Since for any finite dimensional representation $V$,
$\pi_V(c)=0$. It follows from (\ref{normalized-ybe2}) that
$R(z)$ obeys the graded YBE (\ref{rrr}). 
Following Frenkel-Reshetikhin
\cite{Fre92}, we define the `right' dual module $V^*$ and
`left' dual module ${}^*V$ of $V$ by
\beq
\pi_{V^*}(a)=\pi_V(S(a))^{st},~~~~~~
  \pi_{{}^*V}(a)=\pi_V(S^{-1}(a))^{st},\label{dual}
\eeq
respectively. Here $st$ is the supertransposition operation defined by
\beq
(A_{ab})^{st}=(-1)^{[a]([a]+[b])}A_{ba}.
\eeq
Note that in general
$((A_{ab})^{st})^{st}=(-1)^{[a]+[b]}A_{ab}\neq A_{ab}$. Let $ist$ be the
inverse operation of $st$ such that $((A_{ab})^{st})^{ist}=
((A_{ab})^{ist})^{st}=A_{ab}$. Then
\beq
(A_{ab})^{ist}=(-1)^{[b]([a]+[b])}A_{ba}=(-1)^{[a]+[b]}(A_{ab})^{st},
   \label{ist-st}
\eeq
or $A^{ist}=\eta A^{st} \eta$, where $\eta$ is a diagonal matrix with
elements $\eta_{ab}=(-1)^{[a]}\d_{ab}$.

By means of (\ref{dual}) and (\ref{s-r}), one can show that 
\beq
R^{V^*,V}(z)=(R(z)^{-1})^{st_1},~~~~~
  R^{V,{}^*V}(z)=(R(z)^{-1})^{st_2}.
\eeq
{}From the representations for $R^{V^{**},V}(z)$ and
$R^{V,{}^{**}V}(z)$, and the formulae for the action of the square of
antipode:
\beq
S^2(a)=q^{-2h_\rho}\,D_{q^{-2g}}\,(a)\,q^{2h_\rho},~~~~~
       S^{-2}(a)=q^{2h_\rho}\,D_{q^{2g}}\,(a)\,q^{-2h_\rho},~~~~
       \forall a\in U_q[{\cal G}^{(k)}],
\eeq
which can be checked on the generators [remembering that the simple
roots associated with $e_0,~f_0$ are $\a_0=\pm(\frac{1}{k}\d-\psi)$,
respectively, where $\d=(0,0,1)$], it follows that for
any finite dimensional representation $V$,
the R-matrix satisfies the following crossing-unitarity relations:
\bea
(((R(z)^{-1})^{st_1})^{-1})^{st_1}&=&(\pi_V(q^{-2h_\rho})\otimes 1)
    ((R(zq^{-2g}))^{st_1})^{st_1}(\pi_V(q^{2h_\rho})
    \otimes 1),\no\\
(((R(z)^{-1})^{st_2})^{-1})^{st_2}&=&(1\otimes \pi_V(q^{2h_\rho}))
    ((R(zq^{2g}))^{st_2})^{st_2}(1\otimes\pi_V(q^{-2h_\rho})).
    \label{cu-twisted}
\eea
Note also that
\beq
(\pi_V(q^{\pm 2h_\rho})\otimes\pi_V(q^{\pm 2h_\rho}))R(z)
  =R(z)(\pi_V(q^{\pm 2h_\rho})\otimes\pi_V(q^{\pm 2h_\rho})).\label{rho}
\eeq
We introduce the graded permutation operator $P$ on the tensor product
module $V\otimes V$:  $P(v_\a\otimes v_\b)=(-1)^{[\a][\b]}
(v_\b\otimes v_\a)\,,~\forall v_\a, v_\b\in V$. Using the irreducibility
of $V(z)\otimes V$, where $\pi_{V(z)}(a)=\pi_V(D_z(a)) ~\forall a\in
U_q[\Gk]$, and the crossing-unitarity relations (\ref{cu-twisted}),
one can show that the R-matrix $R(z)$ is unitary, that is, 
\beq
R_{12}({z\over w})R_{21}({w\over z})=1,\label{unitarity}
\eeq
where $R_{21}(z)=P_{12}R_{12}(z)P_{12}$.

Now we are going to realize the super RS algebra in $U_q[\Gk]$. Let
\bea
L^+(z)&=&(\pi_V\otimes 1)\R(zq^{\frac{c_2}{2}}),\no\\
L^-(z)&=&(\pi_V\otimes 1)
   \R_{21}(z^{-1}q^{-\frac{c_2}{2}})^{-1},\label{def-L}
\eea
where $c_2=1\otimes c$. It is worth noting that our definition (\ref{def-L})
of $L^\pm(z)$ is different from
that of Frenkel and Reshetikhin \cite{Fre92}. Then, from
(\ref{normalized-ybe2}) and (\ref{unitarity}) we obtain:
\bea
R({z\over w})L_1^\pm(z)L_2^\pm(w)&=&L_2^\pm(w)L_1^\pm(z)R({z\over
         w}),\no\\
R({z_+\over w_-})L_1^+(z)L_2^-(w)&=&L_2^-(w)L_1^+(z)R({z_-\over
         w_+}),\label{super-rs-realization}
\eea
which are nothing but the defining relations (\ref{super-rs})
of the super RS algebra $U(\R)$.
Moreover, from the formulae for the action of
the coproduct and antipode on $\R(z)$ we derive:
\bea
&&\D^T(L^\pm(z)=L^\pm(zq^{\pm 1\otimes {c\over 2}})\stackrel{.}{\otimes}
     L^\pm(zq^{\mp {c\over 2}\otimes 1}),\no\\
&&S^{-1}(L^\pm(z))=L^\pm(z)^{-1}.
\eea
We thus arrive at
\begin{Proposition}: Equations (\ref{def-L}) give a realization of 
$U(\R)$ in $U_q[\Gk]$, with the opposite Hopf algebra structure on
$U_q[\Gk]$.
\end{Proposition}

Now let $R(z)$ be the R-matrix associated with the minimal
$M$-dimensional defining representation $V$ of $U_q(\G)$. We have
\begin{Theorem}\label{df}: $L^\pm(z)$ has the following unique
Gauss decomposition
\bea
L^\pm(z)&=&\left (
\begin{array}{cccc}
1 & \cdots & {} & 0\\
e^\pm_{2,1}(z) & \ddots & {} & {}\\
e^\pm_{3,1}(z) & {}     & {} & \vdots\\
\vdots &  {} & {} & {}\\
e^\pm_{M,1}(z) & \cdots & e^\pm_{M,M-1}(z) & 1
\end{array}
\right )
\left (
\begin{array}{ccc}
k^\pm_1(z) & \cdots & 0\\
\vdots & \ddots & \vdots\\
0 & \cdots & k^\pm_{M}(z)
\end{array}
\right )\no\\
& & \times \left (
\begin{array}{ccccc}
1 & f^\pm_{1,2}(z) & f^\pm_{1,3}(z) & \cdots & f^\pm_{1,M}(z)\\
\vdots & \ddots & \cdots & {} & \vdots\\
{} & {} & {} & {} & f^\pm_{M-1,M}(z)\\
0 & {} & {} & {} & 1
\end{array}
\right ),
\eea
where $e^\pm_{i,j}(z),~f^\pm_{j,i}(z)$ and $k^\pm_i(z) ~(i>j)$ are 
elements in the super RS algebra and $k^\pm_i(z)$ are invertible. 
Let \footnote{Note that our notations here for $X^\pm_i(z)$ are
different from those used
in previous papers \cite{Zha97a,Zha97b} where $X^\pm_i(z)$ are defined
as $X^+_i(z)=e^+_{i+1,i}(z_-)-e^-_{i+1,i}(z_+)$ and 
   $X^-_i(z)=f^+_{i,i+1}(z_+)-f^-_{i,i+1}(z_-)$, respectively.}
\bea
X^+_i(z)&=&f^+_{i,i+1}(z_+)-f^-_{i,i+1}(z_-),\no\\
X^-_i(z)&=&e^-_{i+1,i}(z_+)-e^+_{i+1,i}(z_-),
\eea
where $z_\pm=zq^{\pm{c\over 2}}$, then the defining relations of
$U_q[\Gk]$ can be derived, through 
difference combinations, from relations satisfied by $q^{\pm{c\over 2}}$, $\;
X^\pm_i(z)$, $\; k^\pm_j(z),~ i = 1, 2, \cdots, M-1,\; j = 1, 2, \cdots, M$.
\end{Theorem}

In particular, for the R-matrix $R(z)$ associated with the
$(m+n)$-dimensional representation of $U_q[gl(m|n)]$, we have 
\begin{Theorem} {\rm \cite{Zha97a}}: 
Then the relations satisfied by $\{q^{\pm {c\over 2}}$, $\;
X^\pm_i(z)$, $\; k^\pm_j(z),~i = 1, 2, \cdots, m+n-1,
\;j = 1, 2, \cdots, m+n\}$  are nothing but
the defining relations of $\uqgh$.
\end{Theorem}

\noindent {\bf Remark}:
The theorems are supersymmetric generalizations of that of Ding-Frenkel
\cite{Din93} for the bosonic case.
The Gauss decomposition implies that the elements $e^\pm_{i,j}(z),\;
f^\pm_{j,i}(z) ~(i>j)$ and $k^\pm_i(z)$ are uniquely determined by
$L^\pm(z)$. In the following we will denote $f^\pm_{i,i+1}(z)$,
$e^\pm_{i+1,i}(z)$ as $f^\pm_i(z),\;e^\pm_i(z)$, respectively.

\sect{Ungrading Multiplication Rule of Tensor Products}

We define the matrix elements of $R(z)$ and $L^\pm(z)$ by
\bea
&&R(z)(v_{\a'}\otimes v_{\b'})=R(z)_{\a\b,\a'\b'}(v_\a\otimes v_\b),\no\\
&&L^\pm(z)v_{\a'}=L^\pm(z)_{\a\a'}v_\a.
\eea

In matrix form, (\ref{super-rs}) carries extra signs due to the graded
multiplication rule of tensor products:
\bea
&&R({z\over w})_{\a\b,\a''\b''}L^\pm(z)_{\a''\a'}L^\pm(w)_{\b''
      \b'}\,(-1)^{[\a']([\b']+[\b''])}\no\\
&&~~~~~~~~~~~~~~~=L^\pm(w)_{\b\b''}L^\pm(z)_{\a\a''}R({z\over w})
      _{\a''\b'',\a'\b'}\,(-1)^{[\a]([\b]+[\b''])},\no\\
&&R({z_+\over w_-})_{\a\b,\a''\b''}L^+(z)_{\a''\a'}L^-(w)_{\b''
      \b'}\,(-1)^{[\a']([\b']+[\b''])}\no\\
&&~~~~~~~~~~~~~~~=L^-(w)_{\b\b''}L^+(z)_{\a\a''}R({z_-\over w_+})
      _{\a''\b'',\a'\b'}\,(-1)^{[\a]([\b]+[\b''])}.\label{rll-component}
\eea

We introduce matrix $\t$:
\beq
\t_{\a\b,\a'\b'}=(-1)^{[\a][\b]}\d_{\a\a'}\d_{\b\b'}.\label{theta}
\eeq
With the help of this matrix $\t$, one can cast (\ref{rll-component}) into the
usual matrix equations,
\bea
R({z\over w})L_1^\pm(z)\t L_2^\pm(w)\t&=&\t L_2^\pm(w)\t L_1^\pm(z)R({z\over
         w}),\no\\
R({z_+\over w_-})L_1^+(z)\t L_2^-(w)\t &=&\t L_2^-(w)\t L_1^+(z)R({z_-\over
         w_+}).\label{RLL-LLR1}
\eea
Now the multiplications in ({\ref{RLL-LLR1}) are simply the usual matrix
multiplications, that is the tensor products in (\ref{RLL-LLR1})
carry no grading.

The following matrix equations can be deduced from (\ref{RLL-LLR1}):
\bea
R_{21}({z\over w})\t L^\pm_2(z)\t L_1^\pm(w)&=&
  L^\pm_1(w)\t L^\pm_2(z)\t R_{21}({z\over w}),\label{RLL-LLR2}\\
R_{21}({z_+\over w_-})\t L^+_2(z)\t L_1^-(w)&=&
  L^-_1(w)\t L^+_2(z)\t R_{21}({z_-\over w_+}),\label{RLL-LLR3}\\
R_{21}({z_-\over w_+})\t L^-_2(z)\t L_1^+(w)&=&
  L^+_1(w)\t L^-_2(z)\t R_{21}({z_+\over w_-}),\label{RLL-LLR4}\\
\t L^\pm_2(z)^{-1}\t L^\pm_1(w)^{-1}R_{21}({z\over w})&=&
  R_{21}({z\over w})L^\pm_1(w)^{-1}\t L^\pm_2(z)^{-1}
  \t,\label{RLL-LLR5}\\
\t L^+_2(z)^{-1}\t L^-_1(w)^{-1}R_{21}({z_+\over w_-})&=&
  R_{21}({z_-\over w_+})L^-_1(w)^{-1}\t L^+_2(z)^{-1}
  \t,\label{RLL-LLR6}\\
\t L^-_2(z)^{-1}\t L^+_1(w)^{-1}R_{21}({z_-\over w_+})&=&
  R_{21}({z_+\over w_-})L^+_1(w)^{-1}\t L^-_2(z)^{-1}
  \t,\label{RLL-LLR7}\\
L^\pm_1(w)^{-1}R_{21}({z\over w})\t L^\pm_2(z)\t&=&
  \t L^\pm_2(z)\t R_{21}({z\over
  w})L^\pm_1(w)^{-1},\label{RLL-LLR8}\\
L^-_1(w)^{-1}R_{21}({z_+\over w_-})\t L^+_2(z)\t&=&
  \t L^+_2(z)\t R_{21}({z_-\over
  w_+})L^-_1(w)^{-1},\label{RLL-LLR9}\\
L^+_1(w)^{-1}R_{21}({z_-\over w_+})\t L^-_2(z)\t&=&
  \t L^-_2(z)\t R_{21}({z_+\over
  w_-})L^+_1(w)^{-1}.\label{RLL-LLR10}
\eea
As in (\ref{RLL-LLR1}), the
multiplications in (\ref{RLL-LLR2} -- \ref{RLL-LLR10}) are usual matrix
multiplications.

\sect{Drinfeld Current Realization of $\uqoh$}

We take $R({z\over w})\in End(V\otimes V)$ to be the 
R-matrix with $V$ being the 3-dimensional vector representation of
$U_q(osp(1|2)]$.
Let basis vectors $v_1,\;v_3$ be even and $v_2$ odd.
The R-matrix has the following form:
\beq
R({z\over w})=\left(
\begin{array}{ccccccccc}
1 & 0 & 0 & 0 & 0 & 0 & 0 & 0 & 0\\
0 & a & 0 & b & 0 & 0 & 0 & 0 & 0\\
0 & 0 & d & 0 & c & 0 & r & 0 & 0\\
0 & f & 0 & a & 0 & 0 & 0 & 0 & 0\\
0 & 0 & g & 0 & e & 0 & c & 0 & 0\\
0 & 0 & 0 & 0 & 0 & a & 0 & b & 0\\
0 & 0 & s & 0 & g & 0 & d & 0 & 0\\
0 & 0 & 0 & 0 & 0 & f & 0 & a & 0\\
0 & 0 & 0 & 0 & 0 & 0 & 0 & 0 & 1
\end{array}
\right), \label{r12}
\eeq
where 
\bea
&&a=\frac{q(z-w)}{zq^2-w},~~~~b=\frac{w(q^2-1)}{zq^2-w},~~~~
  c=\frac{q^{1/2}w(q^2-1)(z-w)}{(zq^2-w)(zq^3-w)},\no\\
&&d=\frac{q^2(z-w)(zq-w)}{(zq^2-w)(zq^3-w)},~~~~
  e=a-\frac{zw(q^2-1)(q^3-1)}{(zq^2-w)(zq^3-w)},\no\\
&&f=\frac{z(q^2-1)}{zq^2-w},~~~~
  g=-\frac{q^{5/2}z(q^2-1)(z-w)}{(zq^2-w)(zq^3-w)},\no\\
&&r=\frac{w(q^2-1)[q^3z+q(z-w)-w]}{(zq^2-w)(zq^3-w)},~~~~
  s=\frac{z(q^2-1)[q^3z+q^2(z-w)-w]}{(zq^2-w)(zq^3-w)}.
\eea
$R_{21}(\frac{z}{w})=R(\frac{w}{z})^{-1}$ takes the form
\beq
R_{21}({z\over w})=\left(
\begin{array}{ccccccccc}
1 & 0 & 0 & 0 & 0 & 0 & 0 & 0 & 0\\
0 & a & 0 & f & 0 & 0 & 0 & 0 & 0\\
0 & 0 & d & 0 & -g & 0 & s & 0 & 0\\
0 & b & 0 & a & 0 & 0 & 0 & 0 & 0\\
0 & 0 & -c & 0 & e & 0 & -g & 0 & 0\\
0 & 0 & 0 & 0 & 0 & a & 0 & f & 0\\
0 & 0 & r & 0 & -c & 0 & d & 0 & 0\\
0 & 0 & 0 & 0 & 0 & b & 0 & a & 0\\
0 & 0 & 0 & 0 & 0 & 0 & 0 & 0 & 1
\end{array}
\right). \label{r21}
\eeq
We will construct Drinfeld
current realization of $\uqoh$. We first note that in the present case,
the theorem \ref{df} implies the following decomposition for $L^\pm(z)$:
\bea
L^\pm(z)&=&\left (
\begin{array}{ccc}
1 & 0 &  0\\
e^\pm_1(z) & 1 & 0\\
e^\pm_{3,1}(z) & e^\pm_2(z) & 1
\end{array}
\right )
\left (
\begin{array}{ccc}
k^\pm_1(z) & 0 & 0\\
0 & k^\pm_2(z) & 0\\
0 & 0 & k^\pm_3(z)
\end{array}
\right )
\left (
\begin{array}{ccc}
1 & f^\pm_1(z) & f^\pm_{1,3}(z) \\
0 & 1 &  f^\pm_2(z)\\
0 & 0 & 1
\end{array}
\right )\no\\
&=&\lt(
\begin{array}{ccc}
k^\pm_1(z) & k^\pm_1(z)f^\pm_1(z) & k^\pm_1(z)f^\pm_{1,3}(z)\\
e^\pm_1(z)k^\pm_1(z) & k^\pm_2(z)+e^\pm_1(z)k^\pm_1(z)f^\pm_1(z) & 
    u^\pm\\
e^\pm_{3,1}(z)k^\pm_1(z) & v^\pm & x^\pm
\end{array}
\rt),\label{l+-}
\eea
where,
\bea
u^\pm&=& k^\pm_2(z)f^\pm_2(z)+e^\pm_1(z)k^\pm_1(z)f^\pm_{1,3}(z),\no\\
v^\pm&=&e^\pm_2(z)k^\pm_2(z)+e^\pm_{3,1}(z)
    k^\pm_1(z)f^\pm_1(z),\no\\
x^\pm&=&k^\pm_3(z)+e^\pm_2(z)k^\pm_2(z)f^\pm_2(z)+e^\pm_{3,1}(z)
   k^\pm_1(z)f^\pm_{1,3}(z).
\eea
The inversions $L^\pm(z)^{-1}$ of (\ref{l+-}) are easily seen to be
\beq
L^\pm(z)^{-1}=\lt(
\begin{array}{ccc}
y^\pm & \tilde{x}^\pm & \tilde{u}^\pm\\
\tilde{y}^\pm & k^\pm_2(z)^{-1}+f^\pm_2(z)k^\pm_3(z)^{-1}e^\pm_2(z) &
   -f^\pm_2(z)k^\pm_3(z)^{-1}\\
\tilde{v}^\pm & -k^\pm_3(z)^{-1}
   e^\pm_2(z) & k^\pm_3(z)^{-1}
\end{array}
\rt),
\eeq
where
\bea
\tilde{u}^\pm&=&[f^\pm_1(z)f^\pm_2(z)-f^\pm_{1,3}(z)]k^\pm_3(z)
   ^{-1},\no\\
\tilde{v}^\pm&=&k^\pm_3(z)^{-1}[e^\pm_2(z)e^\pm_1(z)-e^\pm_{3,1}(z)],\no\\
y^\pm&=&k^\pm_1(z)^{-1}+f^\pm_1(z)k^\pm_2(z)^{-1}e^\pm_1(z)+[f^\pm_1(z)
  f^\pm_2(z)-f^\pm_{1,3}(z)]\no\\
& &\times  k^\pm_3(z)^{-1}[e^\pm_2(z)e^\pm_1(z)-
  e^\pm_{3,1}(z)],\no\\
\tilde{x}^\pm&=&-f^\pm_1(z)k^\pm_2(z)^{-1}+[f^\pm_{1,3}(z)-f^\pm_1(z)
  f^\pm_2(z)]k^\pm_3(z)^{-1}e^\pm_2(z),\no\\
\tilde{y}^\pm&=&-k^\pm_2(z)^{-1}e^\pm_1(z)+f^\pm_2(z)k^\pm_3(z)^{-1}
  [e^\pm_{3,1}(z)-e^\pm_2(z)e^\pm_1(z)].
\eea

By means of (\ref{RLL-LLR1}, \ref{RLL-LLR2} -- \ref{RLL-LLR10}) and
after tedious calculations, we derive
\bea
k^\pm_1(z)k^\pm_1(w)&=&k^\pm_1(w)k^\pm_1(z),\no\\
k^+_1(z)k^-_1(w)&=&k^-_1(w)k^+_1(z),\no\\
k^\pm_2(z)k^\pm_2(w)&=&k^\pm_2(w)k^\pm_2(z),\no\\
k^\pm_3(z)k^\pm_3(w)&=&k^\pm_3(w)k^\pm_3(z),\no\\
k^+_3(z)k^-_3(w)&=&k^-_3(w)k^+_3(z),\no\\
k^\pm_1(z)k^\pm_2(w)&=&k^\pm_2(w)k^\pm_1(z),\no\\
\frac{z_\pm-w_\mp}{z_\pm q^2-w_\mp} k^\pm_1(z)k^\mp_2(w)&=&
    \frac{z_\mp-w_\pm}{z_\mp q^2-w_\pm}
    k^\mp_2(w)k^\pm_1(z),\no\\
k^\pm_1(z)k^\pm_3(w)^{-1}&=&k^\pm_3(w)^{-1}k^\pm_1(z),\no\\
\frac{(z_\mp-w_\pm)(z_\mp q-w_\pm)}{(z_\mp q^2-w_\pm)(z_\mp q^3-w_\pm)}
   k^\pm_1(z)k^\mp_3(w)^{-1}&=&
\frac{(z_\pm-w_\mp)(z_\pm q-w_\mp)}{(z_\pm q^2-w_\mp)(z_\pm q^3-w_\mp)}
   k^\mp_3(w)^{-1}k^\pm_1(z),\no\\
\frac{z_\pm-w_\mp q}{z_\pm q-w_\mp}k^\pm_2(z)k^\mp_2(w)&=&
  \frac{z_\mp-w_\pm q}{z_\mp q-w_\pm}k^\mp_2(w)k^\pm_2(z),\no\\
k^\pm_2(z)^{-1}k^\pm_3(w)^{-1}&=&k^\pm_3(w)^{-1}k^\pm_2(z)^{-1},\no\\
\frac{z_\pm-w_\mp}{z_\pm q^2-w_\mp} k^\pm_2(z)^{-1}k^\mp_3(w)^{-1}&=&
    \frac{z_\mp-w_\pm}{z_\mp q^2-w_\pm}
    k^\mp_3(w)^{-1}k^\pm_2(z)^{-1},\label{k1k2k3}
\eea
\bea
k^\pm_1(z)X^-_1(w)k^\pm_1(z)^{-1}&=&\frac{z_\pm q^2-w}{q(z_\pm-w)}
    X^-_1(w),\no\\
k^\pm_1(z)^{-1}X^+_1(w)k^\pm_1(z)&=&\frac{z_\mp q^2-w}{q(z_\mp-w)}
    X^+_1(w),\no\\
k^\pm_2(z)X^-_1(w)k^\pm_2(z)^{-1}&=&\frac{(z_\pm-w q^2)(z_\pm q-w)}{q(z_\pm-w)
    (z_\pm-wq)}X^-_1(w),\no\\
k^\pm_2(z)^{-1}X^+_1(w)k^\pm_2(z)&=&\frac{(z_\mp-w q^2)(z_\mp q-w)}{q(z_\mp-w)
    (z_\mp-wq)}X^+_1(w),\no\\
k^\pm_3(z)X^-_1(w)k^\pm_3(z)^{-1}&=&\frac{z_\pm-wq^3}{q(z_\pm-wq)}
    X^-_1(w),\no\\
k^\pm_3(z)^{-1}X^+_1(w)k^\pm_3(z)&=&\frac{z_\mp-wq^3}{q(z_\mp-wq)}
    X^+_1(w),\no\\
k^\pm_1(z)X^-_2(w)k^\pm_1(z)^{-1}&=&\frac{z_\pm q^3-w}{q(z_\pm q-w)}
    X^-_2(w),\no\\
k^\pm_1(z)^{-1}X^+_2(w)k^\pm_1(z)&=&\frac{z_\mp q^3-w}{q(z_\mp q-w)}
    X^+_2(w),\no\\
k^\pm_2(z)X^-_2(w)k^\pm_2(z)^{-1}&=&\frac{(z_\pm-w q)(z_\pm q^2-w)}
    {q(z_\pm q-w)(z_\pm-w)}X^-_2(w),\no\\
k^\pm_2(z)^{-1}X^+_2(w)k^\pm_2(z)&=&\frac{(z_\mp-w q)(z_\mp q^2-w)}
    {q(z_\mp q-w)(z_\mp-w)}X^+_2(w),\no\\
k^\pm_3(z)X^-_2(w)k^\pm_3(z)^{-1}&=&\frac{z_\pm-wq^2}{q(z_\pm-w)}
    X^-_2(w),\no\\
k^\pm_3(z)^{-1}X^+_2(w)k^\pm_3(z)&=&\frac{z_\mp-wq^2}{q(z_\mp-w)}
    X^+_2(w),\label{x1k1k3-x2k1k3}
\eea
\bea
\frac{z-w}{zq^2-w}X^-_1(z)X^-_2(w)+\frac{z-wq}{zq^3-w}X^-_2(w)X^-_1(z)
   &=&0,\no\\
\frac{z-wq}{zq-w}X^-_1(z)X^-_1(w)+\frac{z-wq^2}{zq^2-w}X^-_1(w)X^-_1(z)
   &=&0,\no\\
\frac{z-wq}{zq-w}X^-_2(z)X^-_2(w)+\frac{z-wq^2}{zq^2-w}X^-_2(w)X^-_2(z)
   &=&0,\no\\
\frac{z-wq}{zq^3-w}X^+_1(z)X^+_2(w)+\frac{z-w}{zq^2-w}X^+_2(w)X^+_1(z)
   &=&0,\no\\
\frac{z-wq^2}{zq^2-w}X^+_1(z)X^+_1(w)+\frac{z-wq}{zq-w}X^+_1(w)X^+_1(z)
   &=&0,\no\\
\frac{z-wq^2}{zq^2-w}X^+_2(z)X^+_2(w)+\frac{z-wq}{zq-w}X^+_2(w)X^+_2(z)
   &=&0,\label{x+x+x-x-}
\eea
\bea
\{X^-_1(w), X^+_1(z)\}&=&(q-q^{-1})\lt[-\d(\frac{z}{w}q^c)k^+_2(z_+)
    k^+_1(z_+)^{-1}\rt.\no\\
& &~~~~~~~~\lt.+\d(\frac{z}{w}q^{-c})k^-_2(w_+)k^-_1(w_+)^{-1}\rt],\no\\
\{X^-_2(w), X^+_2(z)\}&=&(q-q^{-1})\lt[\d(\frac{z}{w}q^c)k^+_3(z_+)
    k^+_2(z_+)^{-1}\rt.\no\\
& &~~~~~~~~\lt.-\d(\frac{z}{w}q^{-c})k^-_3(w_+)k^-_2(w_+)^{-1}\rt],\no\\
\{X^-_2(w), X^+_1(z)\}&=&(q-q^{-1})q^{-{1\over 2}}
    \lt[-\d(\frac{z}{w}q^{c+1})k^+_2(z_+)k^+_1(z_+)^{-1}\rt.\no\\
& &~~~~~~~~\lt.+\d(\frac{z}{w}q^{-c+1})k^-_3(w_+)k^-_2(w_+)^{-1}\rt],\no\\
\{X^-_1(z), X^+_2(w)\}&=&(q-q^{-1})q^{1\over 2}\lt[\d(\frac{w}{z}q^{c-1})
    k^+_2(z_-)k^+_1(z_-)^{-1}\rt.\no\\
& &~~~~~~~~\lt.-\d(\frac{w}{z}q^{-c-1})k^-_3(w_-)k^-_2(w_-)^{-1}\rt],
  \label{x1-x2}
\eea
where $\{X,Y\}\equiv XY+YX$ denotes an anti-commutator, and
\beq
\d(z)=\sum_{l\in {\bf Z}}\,z^l
\eeq
is a formal series. $\d(z)$ enjoys the following properties:
\beq
\d(\frac{z}{w})=\d(\frac{w}{z}),~~~~~~
  \d(\frac{z}{w})f(z)=\d(\frac{z}{w})f(w).
\eeq

The last relation in (\ref{x1-x2}) can be recast into the following form
\bea
\{X^-_1(w), X^+_2(z)\}&=&(q-q^{-1})q^{1\over 2}\lt[\d(\frac{z}{w}q^{c-1})
    k^+_2(z_+q^{-1})k^+_1(z_+q^{-1})^{-1}\rt.\no\\
& &~~~~~~~~\lt.-\d(\frac{z}{w}q^{-c-1})k^-_3(w_+q)k^-_2(w_+q)^{-1}\rt].
  \label{x1w-x2z}
\eea

We define the following difference combinations:
\beq
X^\pm(z)=(q-q^{-1})\lt[X^\pm_1(z)+X^\pm_2(zq)\rt].\label{x=x1+x2}
\eeq
Then (\ref{x1k1k3-x2k1k3} -- \ref{x1w-x2z}) can be rewritten as
\bea
k^\pm_1(z)X^-(w)k^\pm_1(z)^{-1}&=&\frac{z_\pm q^2-w}{q(z_\pm-w)}
    X^-(w),\no\\
k^\pm_1(z)^{-1}X^+(w)k^\pm_1(z)&=&\frac{z_\mp q^2-w}{q(z_\mp-w)}
    X^+(w),\no\\
k^\pm_2(z)X^-(w)k^\pm_2(z)^{-1}&=&\frac{(z_\pm-w q^2)(z_\pm q-w)}{q(z_\pm-w)
    (z_\pm-wq)}X^-(w),\no\\
k^\pm_2(z)^{-1}X^+(w)k^\pm_2(z)&=&\frac{(z_\mp-w q^2)(z_\mp q-w)}{q(z_\mp-w)
    (z_\mp-wq)}X^+(w),\no\\
k^\pm_3(z)X^-(w)k^\pm_3(z)^{-1}&=&\frac{z_\pm-wq^3}{q(z_\pm-wq)}
    X^-(w),\no\\
k^\pm_3(z)^{-1}X^+(w)k^\pm_3(z)&=&\frac{z_\mp-wq^3}{q(z_\mp-wq)}
    X^+(w),\label{x+-k1k2k3}
\eea
\bea
\frac{z-wq}{zq-w}X^-(z)X^-(w)+\frac{z-wq^2}{zq^2-w}X^-(w)X^-(z)
   &=&0,\no\\
\frac{z-wq^2}{zq^2-w}X^+(z)X^+(w)+\frac{z-wq}{zq-w}X^+(w)X^+(z)
   &=&0,\label{x++x--}
\eea
\bea
\{X^-(w), X^+(z)\}&=&\frac{-1}{q-q^{-1}}\lt[\d(\frac{z}{w}q^c)\lt(
   (1+q^{-{1\over 2}}-q^{1\over 2}) k^+_2(z_+)
    k^+_1(z_+)^{-1}-k^+_3(z_+q)k^+_2(z_+q)^{-1}\rt)\rt.\no\\
& &\lt.-\d(\frac{z}{w}q^{-c})\lt(k^-_2(w_+)k^-_1(w_+)^{-1}
   -(1+q^{-{1\over 2}}-q^{1\over 2})k^-_3(w_+q)k^-_2(w_+q)^{-1}\rt)\rt].\no\\
   \label{x+x-}
\eea

Further defining
\bea
&&\phi_i(z)=k^+_{i+1}(z)k^+_i(z)^{-1},\no\\
&&\psi_i(z)=k^-_{i+1}(z)k^-_i(z)^{-1},~~i=1,2,\no\\
&&\phi(z)=(1+q^{-{1\over 2}}-q^{1\over 2})\phi_1(z)-\phi_2(zq),\no\\
&&\psi(z)=\psi_1(z)-(1+q^{-{1\over 2}}-q^{1\over 2})\psi_2(zq),\label{phi-psi}
\eea
then we have
\begin{Theorem}:  $q^{\pm{c\over 2}},\;
X^\pm(z),\;\phi(z),\;\psi(z)$
give the defining relations of $\uqoh$. More precisely,
$\uqoh$ is an associative algebra with unit 1 and
the Drinfeld generators: $X^\pm(z),~\phi(z)$ and $\psi(z)$, a central
element $c$ and a nonzero complex parameter $q$. $\phi(z)$ and $\psi(z)$
are invertible. The gradings of the generators are: $[X^\pm(z)]=1$ and
$[\phi(z)]=[\psi(z)]=[c]=0$. The relations are given by
\begin{eqnarray*}
\phi(z)\phi(w)&=&\phi(w)\phi(z),\no\\
\psi(z)\psi(w)&=&\psi(w)\psi(z),\no\\
\phi(z)\psi(w)\phi(z)^{-1}\psi(w)^{-1}&=&\frac{(z_+q-w_-)(z_--w_+q)
    (z_+-w_-q^2)(z_-q^2-w_+)}{(z_+-w_-q)(z_-q-w_+)(z_+q^2-w_-)
    (z_--w_+q^2)},
\end{eqnarray*}
\begin{eqnarray*}
\phi(z)X^-(w)\phi(z)^{-1}&=&\frac{(z_+-w q^2)(z_+ q-w)}{(z_+ q^2-w)
    (z_+-wq)}X^-(w),\no\\
\phi(z)^{-1}X^+(w)\phi(z)&=&\frac{(z_--w q^2)(z_- q-w)}{(z_- q^2-w)
    (z_--wq)}X^+(w),\no\\
\psi(z)X^-(w)\psi(z)^{-1}&=&\frac{(z_--w q^2)(z_- q-w)}{(z_- q^2-w)
    (z_--wq)}X^-(w),\no\\
\psi(z)^{-1}X^+(w)\psi(z)&=&\frac{(z_+-w q^2)(z_+ q-w)}{(z_+ q^2-w)
    (z_+-wq)}X^+(w),\label{x-phipsi}
\end{eqnarray*}
\begin{eqnarray*}
\frac{z-wq}{zq-w}X^-(z)X^-(w)+\frac{z-wq^2}{zq^2-w}X^-(w)X^-(z)
   &=&0,\no\\
\frac{z-wq^2}{zq^2-w}X^+(z)X^+(w)+\frac{z-wq}{zq-w}X^+(w)X^+(z)
   &=&0,\label{Fx++x--}
\end{eqnarray*}
\beq
\{X^+(z), X^-(w)\}=\frac{1}{q-q^{-1}}\lt[
     \d(\frac{w}{z}q^{c})\psi(w_+)
     -\d(\frac{w}{z}q^{-c})\phi(z_+)\rt].\label{x+x-=psiphi}
\eeq
\end{Theorem}

Super-Yangian doubles with center $DY_\hbar[\Gk]$
are the degenerated cases of the
corresponding quantum affine superalgebras $U_q[\Gk]$. The Drinfeld
realization of $DY_\hbar[\Gk]$ can also be obtained  by means of the super RS
algebra with a rational R-matrix $R(z)$ and the Gauss decomposition theorem.
In \cite{Zha97b}, we gave the defining relations of
$DY_\hbar[gl(m|n)^{(1)}]$ in terms of Drinfeld current generators (see
also \cite{Cai97Y} for the simplest case of $m=n=1$).
We now introduce $DY_\hbar[osp(1|2)^{(1)}]$. 
\begin{Theorem}: $DY_\hbar[osp(1|2)^{(1)}]$ is an associative algebra 
over the ring of formal power series in the variable $\hbar$ and
the Drinfeld generators: $X^\pm(u),~\phi(u)$ and $\psi(u)$, and a central
element $c$. $\phi(u)$ and $\psi(u)$
are invertible. The gradings of the generators are: $[X^\pm(z)]=1$ and
$[\phi(z)]=[\psi(z)]=[c]=0$. The defining relations are given by
\begin{eqnarray*}
\phi(u)\phi(v)&=&\phi(v)\phi(u),\no\\
\psi(u)\psi(v)&=&\psi(v)\psi(u),\no\\
\phi(u)\psi(v)\phi(u)^{-1}\psi(v)^{-1}&=&\frac{(u_+-v_-+\hbar)(u_--v_+-\hbar)
    (u_+-v_--2\hbar)(u_--v_++2\hbar)}{(u_+-v_--\hbar)(u_--v_++\hbar)
    (u_+-v_-+2\hbar)(u_--v_+-2\hbar)},
\end{eqnarray*}
\begin{eqnarray*}
\phi(u)X^-(v)\phi(u)^{-1}&=&\frac{(u_+-v-2\hbar)(u_+-v+\hbar)}
    {(u_+-v+2\hbar)(u_+-v-\hbar)}X^-(w),\no\\
\phi(u)^{-1}X^+(v)\phi(u)&=&\frac{(u_--v-2\hbar)(u_--v+\hbar)}
    {(u_--v+2\hbar)(u_--v-\hbar)}X^+(w),\no\\
\psi(u)X^-(v)\psi(u)^{-1}&=&\frac{(u_--v-2\hbar)(u_--v+\hbar)}
    {(u_--v+2\hbar)(u_--v-\hbar)}X^-(w),\no\\
\psi(u)^{-1}X^+(v)\psi(u)&=&\frac{(u_+-v-2\hbar)(u_+-v+\hbar)}
    {(u_+-v+2\hbar)(u_+-v-\hbar)}X^+(w),
\end{eqnarray*}
\begin{eqnarray*}
\frac{u-v-\hbar}{u-v+\hbar}X^-(u)X^-(v)+\frac{u-v-2\hbar}{u-v+2\hbar}
   X^-(v)X^-(u)&=&0,\no\\
\frac{u-v-2\hbar}{u-v+2\hbar}X^+(u)X^+(v)+\frac{u-v-\hbar}{u-v+\hbar}
   X^+(v)X^+(u)&=&0,
\end{eqnarray*}
\beq
\{X^+(u), X^-(v)\}=\frac{1}{2\hbar}\lt[
     \d(u_--v_+)\psi(v_+)
     -\d(u_+-v_-)\phi(u_+)\rt],
\eeq
where $u_\pm=u\pm \frac{1}{2}\hbar c$ and
\beq
\d(u-v)=\sum_{l\in{\bf Z}}\;u^lv^{-l-1}
\eeq
is a formal series.
\end{Theorem}

\vskip.3in
\noindent {\bf Acknowledgements.} Y.-Z.Z. is supported by the 
Queen Elizabeth II Fellowship Grant from Australian Research Council.

%\newpage

\end{document}